\documentclass[12pt]{article}
\usepackage{ifthen}
\usepackage[dvips]{graphicx}
\newcommand\FIGEPS[4]{
{\begin{figure}\begin{center}\includegraphics[scale=0.5]{#1.eps}
\caption{\footnotesize#4}\label{#1}\end{center}\end{figure}}
}

\def\bea{\begin{eqnarray}}
\def\eea{\end{eqnarray}}

\def\beq{\begin{equation}}
\def\eeq{\end{equation}}
\def\ba{\beq\new\begin{array}{c}}
\def\ea{\end{array}\eeq}
\def\be{\ba}
\def\ee{\ea}

\makeatletter
\newdimen\normalarrayskip 
\newdimen\minarrayskip 
\normalarrayskip\baselineskip \minarrayskip\jot
\newif\ifold \oldtrue \def\new{\oldfalse}
\def\arraymode{\ifold\relax\else\displaystyle\fi} 
\def\eqnumphantom{\phantom{(\theequation)}} 
\def\@arrayskip{\ifold\baselineskip\z@\lineskip\z@
\else \baselineskip\minarrayskip\lineskip2\minarrayskip\fi}
\def\@arrayclassz{\ifcase \@lastchclass \@acolampacol \or
\@ampacol \or \or \or \@addamp \or \@acolampacol \or
\@firstampfalse \@acol \fi \edef\@preamble{\@preamble \ifcase
\@chnum \hfil$\relax\arraymode\@sharp$\hfil \or
$\relax\arraymode\@sharp$\hfil \or
\hfil$\relax\arraymode\@sharp$\fi}}
\def\@array[#1]#2{\setbox\@arstrutbox=\hbox{\vrule
height\arraystretch \ht\strutbox depth\arraystretch \dp\strutbox
width\z@}\@mkpream{#2}\edef\@preamble{\halign \noexpand\@halignto
\bgroup \tabskip\z@ \@arstrut \@preamble \tabskip\z@ \cr}%
\let\@startpbox\@@startpbox \let\@endpbox\@@endpbox
\if #1t\vtop \else \if#1b\vbox \else \vcenter \fi\fi \bgroup
\let\par\relax
\let\@sharp##\let\protect\relax
\@arrayskip\@preamble}
\def\eqnarray{\stepcounter{equation}%
\let\@currentlabel=\theequation
\global\@eqnswtrue \global\@eqcnt\z@ \tabskip\@centering
\let\\=\@eqncr
$$%
\halign to \displaywidth\bgroup
\eqnumphantom\@eqnsel\hskip\@centering
$\displaystyle \tabskip\z@ {##}$%
\global\@eqcnt\@ne \hskip 2\arraycolsep
$\displaystyle\arraymode{##}$\hfil \global\@eqcnt\tw@ \hskip
2\arraycolsep $\displaystyle\tabskip\z@{##}$\hfil
\tabskip\@centering &{##}\tabskip\z@\cr}
\begingroup\ifx\undefined\newsymbol \else\def\input#1 {\endgroup}\fi

\begin{document}

\setcounter{footnote}{1}
\def\thefootnote{\fnsymbol{footnote}}

\begin{center}
\vspace{0.3in}

{\Large\bf }
\end{center}

\vspace{1.5cm}

\bigskip


\begin{center}
{\large \textbf{OPERATIONS AND IDENTITIES IN TENSOR ALGEBRA}}
\end{center}

\bigskip

\centerline{{\large Chernyakov Yu.} \footnote{ITEP, Moscow,
Russia, e-mail: chernyakov@itep.ru }, {\large Dolotin V.}
\footnote{ITEP, Moscow, Russia, e-mail: vd@itep.ru }}

\bigskip

\vspace{1.5cm}

\centerline{ \bf{Abstract}}

{\footnotesize We define the class of non-decomposable
$N$-ary operations in the mixed tensor algebra
$\bigoplus\limits_{i,j=0}^\infty A_i^j$. There are higher
Jacobi-like identities for (binary) deformed matrix commutator and a
3-ary operation which is non-decomposable into binary ones.}

\section{Introduction}
The subject of the algebra is mostly deals with studying the action of the rings of linear maps on the corresponding modules ([1],[2]). Superposition
of these maps stipulates the binarity of the composition
law in the space of operators. Linear combinations of such superpositions
give new (generically, non-associative) composition laws, such as
multiplications in Lie algebra, Jordan algebra and their
deformations. Let tensor algebra be generated by a given linear
space $V$. From the tensor algebra point of view the
operator ring is the space $V\otimes V^\ast$ of (1,1)-type
tensors. The composition law in such space is generated by
elementary operations for generators:

$$(e_i\otimes e_j^*)
\cdot(e_k\otimes e_l^*)=(e^*_j,e_k)e_i\otimes
e_l^*=\delta_{jk}\,e_i\otimes e_l^*$$

However, the class of such elementary operations in the mixed
tensor algebra $A^\bullet_\bullet:=\bigoplus\limits_{i,j=0}^\infty
A_i^j=\bigoplus\limits_{i,j=0}^\infty V^{\otimes i}\otimes
V^{*\otimes j}$ is considerably wider. Combinations of these
elementary operations generate the class of operations which in
general is non-binary and is not decomposable into binary
operations. So, there appears a problem of finding such subspaces
in $A^\bullet_\bullet$, which are invariant with respect to these
operations, and finding the identities which these operations
satisfy.

This paper is organized in the following way. In Chapter 2 we
describe the language which simplifies the further consideration. In
Chapter 3 the finite dimensional subalgebras of mixed tensor algebra
$\bigoplus\limits_{i,j=0}^\infty A_i^j$, i.e. the (1,1)-type tensors
(matrices), binary operations and the corresponding identities are
considered. In Chapter 4 we consider the tensor space $A^{1}_{2}
\oplus A^{2}_{1}$, a 3-ary operation non-decomposable into binary operations
and an identity which generalizes the commutator and Jacobi identity. In Chapter 5 we speak about the $n$-ary operations. The Chapter 6 is a discussion.

\section{Language}
Let us consider the action of the set of linear operators on the
vector space $V$. The dual space $V^{\ast}$ consists of linear
forms on $V$. It is possible to define the operation called
contraction, given the value of each element of $V^{\ast}$ on the elements of $V$. The result of this operation is scalar. It is possible to define
another operation called tensor product of $V^{\ast}$ and $V$. The
linear operator is defined as (1,1)-type tensor where digits
denote the number of "upper and lower indices" of the tensor correspondingly, i.e. $A^{1}_{1} = V \otimes V^{\ast}$.
The spaces $V^{\ast}$ and $V$ generate the tensor algebra. We
consider it as a graded space. The value of this
grading is defined as the difference between the number of lower and upper
indexes of the tensor. Thus, the grading of $V$ is $[V] =
+1$ and of $V^{\ast}$ is $[V^{\ast}] = -1$.

It is convenient to describe the general situation by the
following picture: \be

... \oplus

\left[
\begin{array}{c}
 A^{0}_{1} \\
 \oplus \\
 A^{1}_{2} \\
 \oplus \\
 A^{2}_{3} \\
 \oplus \\
 A^{3}_{4} \\
 ... \end{array}
\right]

\oplus

\left[
\begin{array}{c}
 A^{0}_{0} \\
 \oplus \\
 A^{1}_{1} \\
 \oplus \\
 A^{2}_{2} \\
 \oplus \\
 A^{3}_{3} \\
 ... \end{array}
\right]

\oplus

\left[
\begin{array}{c}
 A^{1}_{0} \\
 \oplus \\
 A^{2}_{1} \\
 \oplus \\
 A^{3}_{2} \\
 \oplus \\
 A^{4}_{3} \\
 ... \end{array}
\right]

\oplus ...
 \ee
Let us introduce some diagrammatic language which helps us to
produce further results. We use the notation
$\stackrel{a}{\rightarrow}$ for the (1,1)-type tensor where the
arrow tip means the upper index and the tip without the arrow
means the lower index of the tensor. Using this rule we get for
the tensors $a\in A^{1}_{3}$ and $b\in A^{5}_{2}$ correspondingly
(Fig.1).

\FIGEPS{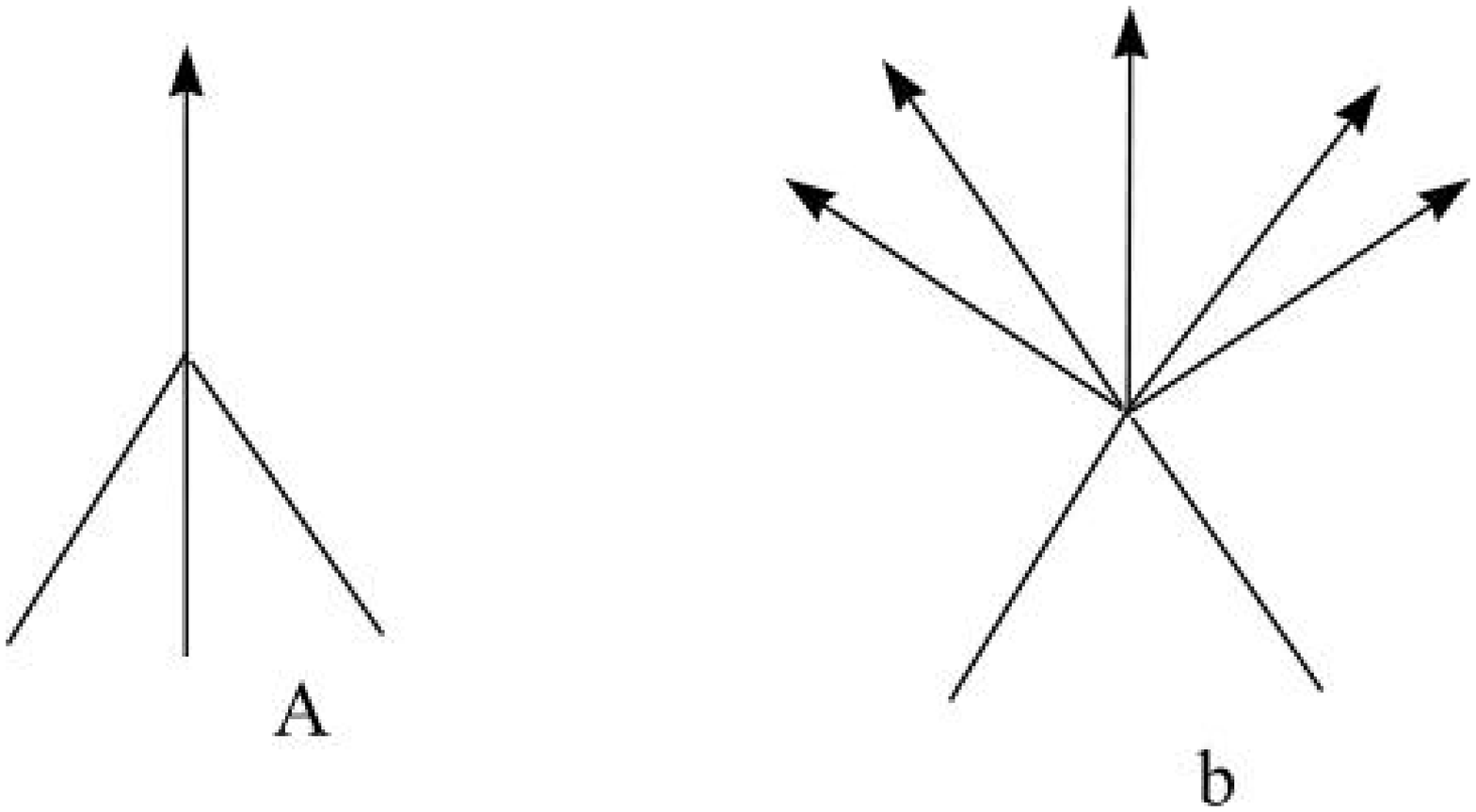}{}{}{}

Following this notation the contraction of two (1,1)-type tensors
$A$ and $B$ using the upper index of $A$ and the lower of $B$
looks as follows figure 2.

As a result we get a (1,1)-type tensor again. Here and further we
take into account that all indices have equal number of values.
So, the contraction between arbitrary two indexes of different
types (upper and lower) is possible. For example, the tensor
product $\otimes$ is represented as figure 3.

\FIGEPS{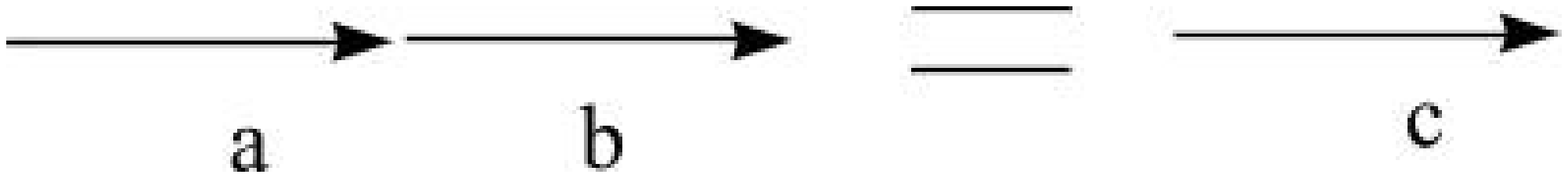}{}{}{}

\FIGEPS{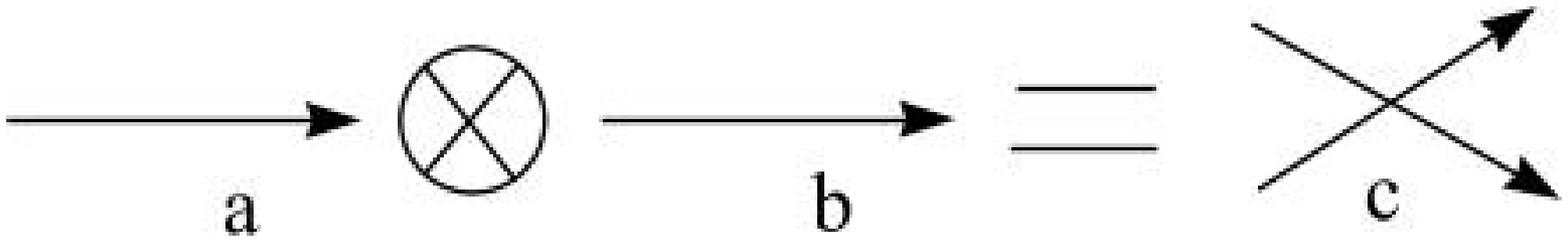}{}{}{}

In this case we get (2,2)-type tensor. Note, that similar notations
were used in the paper [3], but being applied to differential operators, which may be regarded a the subclass of symmetric tensors in $A_\bullet^\bullet$.

In the mixed tensor algebra $A_\bullet^\bullet$
there exist subspaces which are invariant with respect to some fixed
operations (non necessarily binary). These subspaces may be called subalgebras of $A_\bullet^\bullet$. In turn, these operations may satisfy certain identities.

\section{(1,1)-type tensors - matrices and binary operations}

Let us consider the following question: how many different types
of contractions (called further elementary operation) of two
(1,1)-type tensors is there. We consider as the result of the
contraction tensors having the same type as the type of
ingredients and tensors of other types as well. One can check that
using the contraction and tensor product (which corresponds to the
contraction of empty set of indices) we can get seven variants of
difference compositions (Fig.4).

The first four variants are the (1,1)-type tensors, the next two
ones are the (0,0)-type tensors - scalars, and the latter one is
(2,2)-type tensor. Note that the result of each elementary
operation has the grading equal to zero. Evidently this is the
property of contraction only for $(1,1)$-type tensors,
 or other zero graded ($(m,m)$-type) tensors.

Further in this Chapter we will be interested in such elementary
operations having (1,1)-type tensors as the result. These
operations transfer the initial space, the space of (1,1)-type
tensors, to itself. Now let us define a generic binary operation
as a linear composition of elementary operations between tensors
$A$ and $B$ (here we do not write tensor indices applied to
(1,1)-type tensors): \be

A \circ B = \alpha \cdot AB + \ \beta \cdot BA
 +  \gamma \cdot A \ Tr B
 +  \delta \cdot B \ Tr A

\ee where $\alpha, \beta, \gamma, \delta \in \textit{C}$. So, we
define the binary bilinear operation $\varphi_{2}$ on the matrix
space: \be

\varphi_{2}: A_1^1 \times A_1^1 \rightarrow A_1^1

\ee Note, if we put  $\alpha + \beta = 0,\ \gamma =0,\ \delta =0$,
we get \be

A \circ B = \ \alpha \cdot (AB - BA),

\ee

\FIGEPS{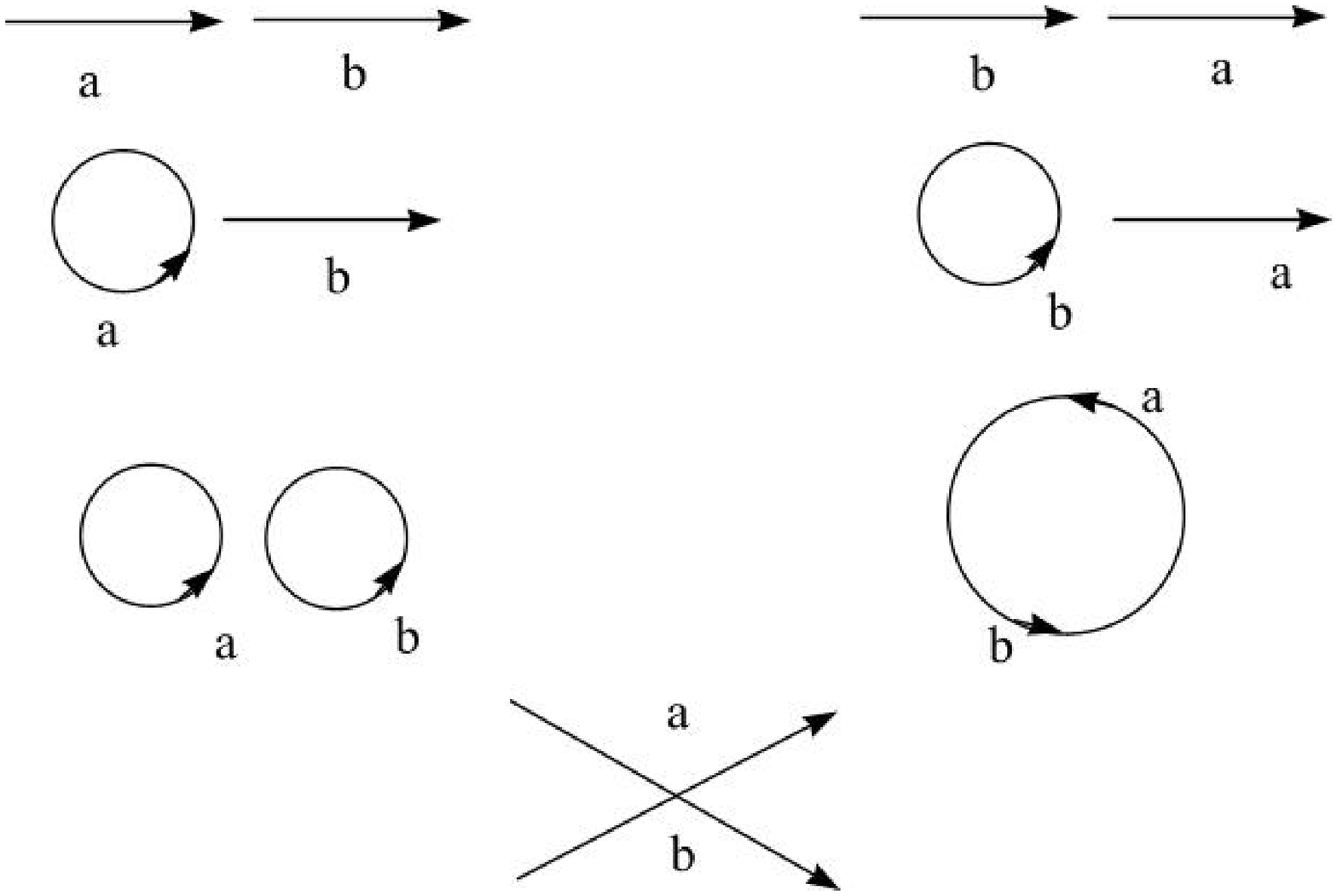}{}{}{}

It is the matrix commutator if we put $\alpha = +1,\ \beta = -1$.
How can we see the Jacobi identity using the language defined in
Chapter 2 ? We have to consider the three-symbol words, for
example $a \ b \ c$ where each symbol corresponds to the
(1,1)-type tensor and each word is the result of contractions of
these tensors (Fig.5).

 \FIGEPS{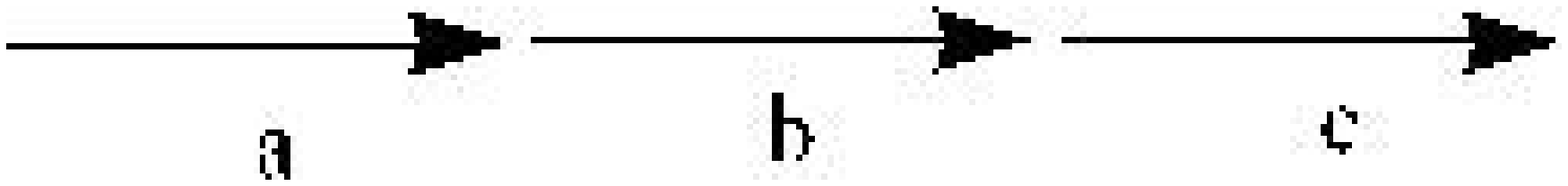}{}{}{}

It is possible to get this word by two ways using the primary
operations: $A \circ B$ and $B \circ C$. On the second step of
this procedure we must operate with the result of the first
primary operation and the rest third tensor. Organizing the
enumeration in the appropriate way (using the permutations) we get
the Jacobi identity: \be

(A \circ B) \circ C + (B \circ C) \circ A + (C \circ A) \circ B =
0.

\ee For the operation (2) it is possible to get the identity
setting $\beta=-\alpha,\ \delta=-\gamma$. This identity is the
generalization
 of the Jacobi identity (see details in Appendix 1 ):
\be

((A \circ B) \circ C) \circ D + ((C \circ B) \circ D) \circ A +
((C \circ D) \circ A) \circ B + ((A \circ D) \circ B) \circ C +\\
((C \circ A) \circ B) \circ D + ((D \circ C) \circ B) \circ A +
((A \circ C) \circ D) \circ B + ((B \circ A) \circ D) \circ C +\\
((B \circ C) \circ A) \circ D + ((B \circ D) \circ C) \circ A +
((D \circ A) \circ C) \circ B + ((D \circ B) \circ A) \circ C \\
= 0.

\ee In order to get this identity we use the following arguments.
First of all, the direct testing shows that ordinary cyclic
permutations of symbols in $(A \circ B) \circ C)$ does not give
zero result, so we can take triple compositions of the operation
(2). In this case we must take into account all primary operations
taking part in formation of four-symbolic words. There are twelve
such primary operations: $(AB), (BA), (CA), (AC), (AD), (DA),
(BC), (CB), (BD), (DB),\\ (DC), (CD)$. It is possible to rewrite
the identity (6) in more adequate and general form. For a
$4$-tuple of elements $A_{1}, A_{2}, A_{3}, A_{4}$ we get
 the following trilinear operation $\varphi_{3}$:
\be

\varphi_{3}(A_{1}, A_{2}, A_{3}) = \\
\sum^{3}_{i=1}(-1)^{i-1} \varphi_{2} ( \varphi_{2} (\overline{A_{i}}), A_{i}) =  \{ \{ A_{2}, A_{3} \}, A_{1} \}
- \{ \{ A_{1}, A_{3} \}, A_{2} \} + \{ \{ A_{1}, A_{2} \}, A_{3},
\}

\ee where $\overline{A_i}$ denotes the complement to $A_i$ in the
triple $A_1,A_2,A_3$, $\varphi_{2}$ is the primary binary
operation (2). Using $\varphi_{2}$ and $\varphi_{3}$, we get: \be

\varphi_{4}(A_{1}, A_{2}, A_{3}, A_{4}) =
\sum^{4}_{i=1}(-1)^{i-1} \varphi_{2}(\varphi_{3}(
\overline{A_i}), A_{i})
\ee
And for $\varphi_2$ with $\alpha + \beta = 0,\ \gamma + \delta =0$ we get
\be
\varphi_{4}(A_{1}, A_{2}, A_{3}, A_{4})\equiv 0.

\ee
 Let us formulate the following

\textbf{Proposition 1}: \textit{ For the operation (2) with
$\alpha + \beta = 0,\ \gamma + \delta =0$ it is possible to define
the combinatoric identity (9) using only the operations (2).}

Note that the identity (9) uses higher operations obtained from
the binary operations (2). This is similar to the existence of
Jacobi identity in matrix algebra defined through the commutation
operation.


\section{The tensor space ${\cal A}^{1}_{2} \oplus {\cal A}^{2}_{1}$
 and 3-ary operation}

Now let us take a pair of tensors $A,B$ in general having a
non-zero gradings $[A]=n, [B] = k$. For any type of contractions
the result is always in the space with grading $n+k$. The
naturally arising question is what can be done in order to remain
in initial grading? The answer is that if we take a space of
linear combination of tensors having different gradings there are
operations with the result being tensor from this original space,
so the grading type is preserved.

Let us make the following observation. For elements $A\in {\cal
A}^{1}_{2} \oplus{\cal  A}^{2}_{1}$ of the direct sum of pure
graded spaces it is possible to find an  operation keeping the
result in initial space. To do this we take two more such linear
combinations $B,C\in {\cal A}^{1}_{2} \oplus{\cal  A}^{2}_{1}$ and
consider various collections of triples of tensors, each being a
pure graded component in the sums $A=A_2^1+A_1^2$,$B=B_2^1+B_1^2$
and $C=C_2^1+C_1^2$, with a set of appropriate contractions. We
exclude the contractions of a tensor with itself. These conditions
define the following rule: in order to get a tensor from the space
${\cal A}^{2}_{1}$ one need to take two (2,1)-type tensors and one
(1,2)-type tensor and for getting a tensor from the space ${\cal
A}^{1}_{2}$ one need to take two (1,2)-type tensors and one
(2,1)-type tensor. The number of elementary operations is
fourteen: seven for each type of resulting tensors (Fig.6 and
Fig.7).

\FIGEPS{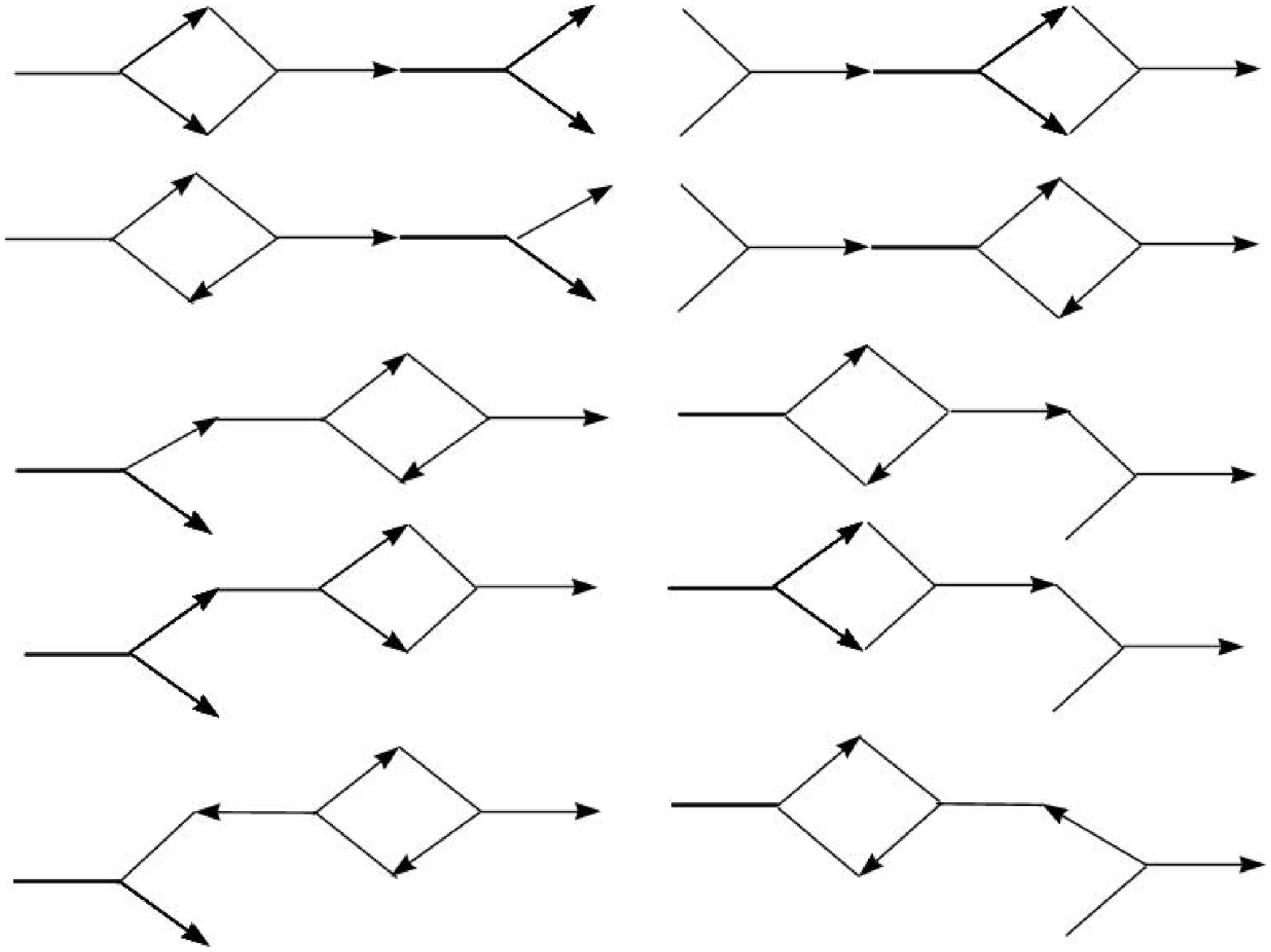}{}{}{}

 We see that in contrast to the matrix
multiplication in the space ${\cal A}^{1}_{1}$ in the case of
tensor space ${\cal A}={\cal A}^{1}_{2} \oplus {\cal A}^{2}_{1}$
the important role plays the type of contractions, in other words
the topology of the corresponding diagram. The diagrams of the
first two elementary operations may be aligned (schematically) in
a line because the objects in this case are placed one after
another. We called these operations by operations of "linear"
type. The diagrams of the other elementary operations can be
placed on the plane.

Note one important property of the operation on the space ${\cal
A}^{1}_{2} \oplus {\cal A}^{2}_{1}$ - they are non-decomposable
into binary operations because the result of any binary operation
may not belongs to the initial space ${\cal A}^{1}_{2} \oplus
{\cal A}^{2}_{1}$. The origin of these 3-ary operations may be
traced from the rule of the action of ${\cal A}^{1}_{2} \oplus
{\cal A}^{2}_{1}$ on the composed "vector space" $V^{\otimes
2}\oplus V^{\otimes 1}={\cal A}^2\oplus{\cal A}^1$, which is
analogous to the linear representation for matrix algebra:
$$(A_2^1+A_1^2)(v^2+v^1)=A_2^1\cdot v^2+A_1^2\cdot v^1$$
So under this action the space $V^{\otimes 2}\oplus V^{\otimes 1}$
passes to itself.

So, taking a linear combination of the elementary operations we
get a generic 3-ary operation: \be
\varphi: {\cal A}\times {\cal A} \times {\cal A} \longrightarrow {\cal A},\\
\ee
while for any binary operation
\be
\tau({\cal A} \times{\cal A})\not\subset{\cal A}.
\ee

\FIGEPS{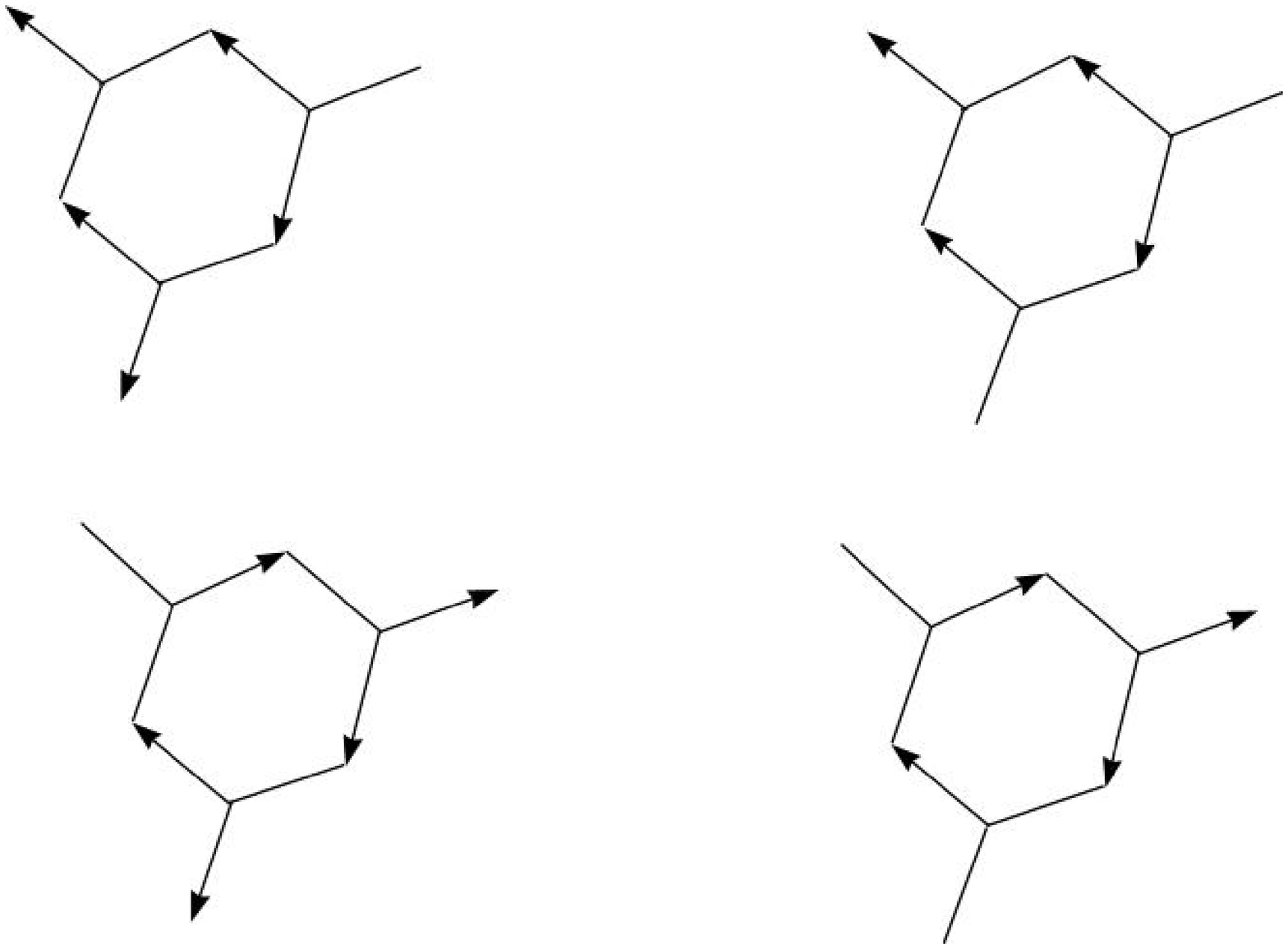}{}{}{}

Now let taking a linear combination of two first elementary
operations we define a 3-ary operation and a Jacobi-like identity.
Let us define the 3-commutator with the condition $\alpha + \beta
+ \gamma = 0$, which is analogous to the condition $\alpha + \beta
= 0$ for the commutator (Fig.8): \be

(X,Y,Z) =\\
 \alpha \cdot X^{2}_{1}Y^{1}_{2}Z^{2}_{1} + \alpha
\cdot Z^{2}_{1}Y^{1}_{2}X^{2}_{1} + \beta \cdot
Z^{2}_{1}X^{1}_{2}Y^{2}_{1} + \beta \cdot
Y^{2}_{1}X^{1}_{2}Z^{2}_{1} + \gamma \cdot
X^{2}_{1}Z^{1}_{2}Y^{2}_{1} + \gamma \cdot
Y^{2}_{1}Z^{1}_{2}X^{2}_{1} +\\
 + \alpha \cdot
X^{1}_{2}Y^{2}_{1}Z^{1}_{2} + \alpha \cdot
Z^{1}_{2}Y^{2}_{1}X^{1}_{2} + \beta \cdot
Z^{1}_{2}X^{2}_{1}Y^{1}_{2} + \beta \cdot
Y^{1}_{2}X^{2}_{1}Z^{1}_{2} + \gamma \cdot
X^{1}_{2}Z^{2}_{1}Y^{1}_{2} + \gamma \cdot
Y^{1}_{2}Z^{2}_{1}X^{1}_{2},

\ee

\FIGEPS{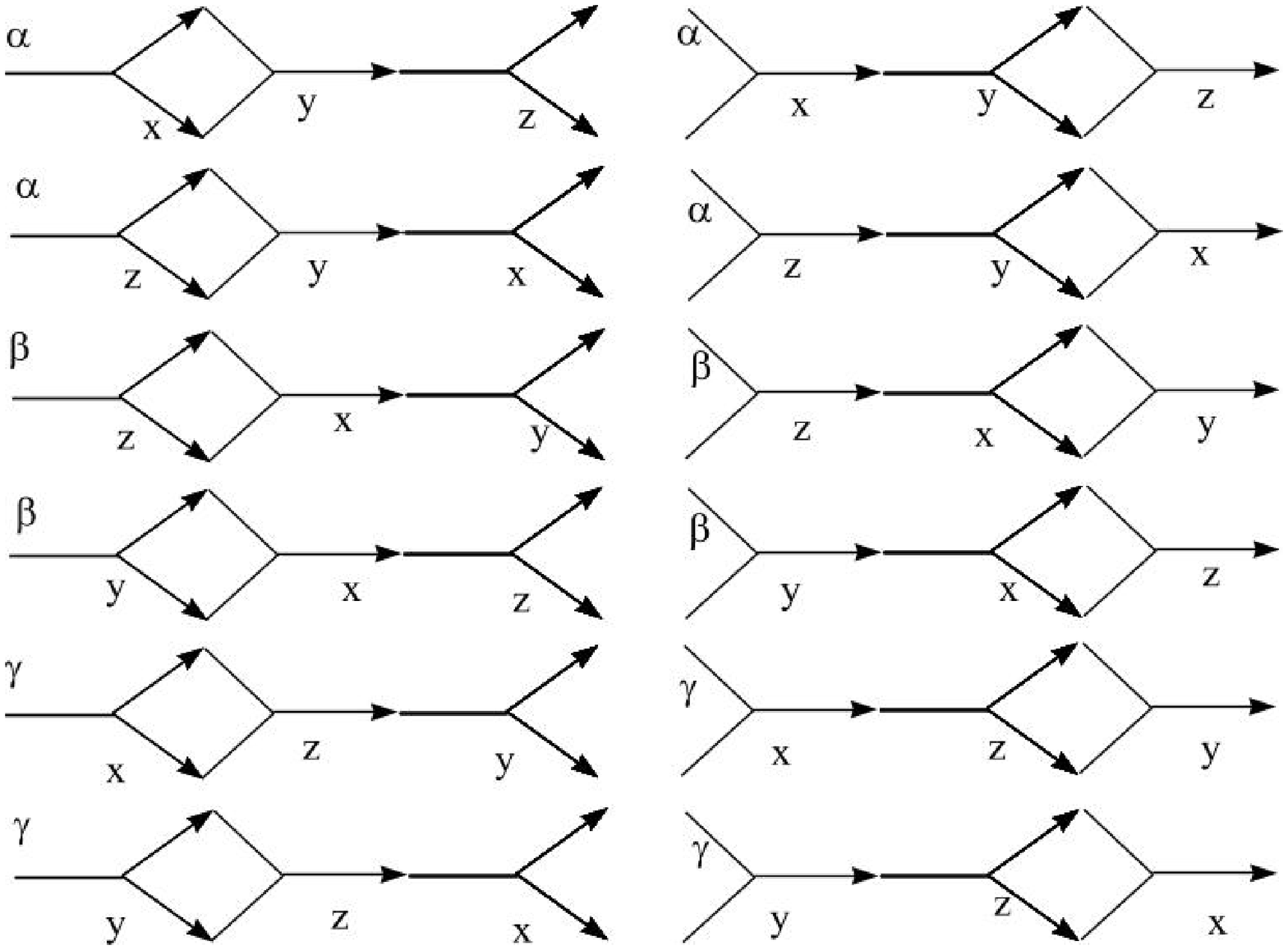}{}{}{}

where: \be

X = X^{1}_{2}+ X^{2}_{1}, \ \ \

Y = Y^{1}_{2}+Y^{2}_{1}, \ \ \

Z = Z^{1}_{2}+Z^{2}_{1}, \ \ \

\ee Here is assumed the summation of the different combinations of
tensor contractions. So, the number of terms in the 3-commutator
is 24. But we can insert the 3-commutator under the sign of
summation and in further we do not specially stipulate the
appeared adding terms from the summing the different tensor
contractions. We set $\alpha, \ \beta, \ \gamma$ equal to the
different values of the cubic root on unity:
\be \alpha =1, \
\beta=e^{i\frac{2\pi}{3}}, \ \gamma=e^{i\frac{4\pi}{3}}. \ee Thus
we can consider 3-commutator modulo any number $\alpha, \ \beta, \
\gamma$. This is analogue of $\alpha = - \beta =1$. By analogy to
the antisymmetry of the commutator:
\be [X,Y] + [Y,X] = 0, \ee the
3-commutator has the following property:
\be (X,Y,Z) + (Z,X,Y) +
(Y,Z,X) = 0, \ee As the result of $n$-times action of the
3-commutator on the tensor space ${\cal A}^{1}_{2} \oplus {\cal
A}^{2}_{1}$ we get the linear combination of the different words
of $2n+1$ symbols. After the second action of the 3-commutator
there appear words consisting of five symbols, like
$A^{2}_{1}B^{1}_{2}C^{2}_{1}D^{1}_{2}E^{2}_{1}$,
$A^{1}_{2}B^{2}_{1}C^{1}_{2}D^{2}_{1}E^{1}_{2}$ and so on. Using
the 3-commutator it is possible to define the identity which
consist of all such words. Let us take as initial the word
$A^{2}_{1}B^{1}_{2}C^{2}_{1}D^{1}_{2}E^{2}_{1}$. Note that this
word can be obtained from the three different primary 3-commutator
operations:
 \be

(A,B,C), \ \ (B,C,D), \ \ (C,D,E).

\ee It is possible to get this operation by the cyclic permutation
in the primary operation $(A,B,C)$, and the secondary one
$((A,B,C)D,E)$. In this case the initial word is
$A^{2}_{1}B^{1}_{2}C^{2}_{1}D^{1}_{2}E^{2}_{1}$. Thus we set the
order of writing the identity although we can take any operation
as a primary one. In a similar way, dealing with the Jacobi
identity we take as initial the word $ABC$ and the operation
$[A,B]$ as primary, and then make cyclic permutations. One needs
to take all primary 3-symbol operations. In other words, all sets
which consist of three symbols chosen from the given five ones. So
we have 10 primary operations modulo the coefficients (weights).
It is possible to get all these operations from the two 5-symbol
classes: $((A,B,C)D,E)$ and $((D,A,C)E,B)$. We must also take into
account the existence of classes which differ by cyclic
permutations in the sets of three symbols:  $(X,Y,Z) + (Z,X,Y) +
(Y,Z,X) = 0$ - $(Z,Y,X) + (X,Z,Y) + (Y,X,Z) = 0$. For a chosen
initial word it does not matter which operation we define as the
primary one: $(A,B,C)$ or $(C,B,A)$. If we start from $(C,B,A)$,
it is possible to get the same sets of tree symbol as for another
four primary operations (as in the case of $(A,B,C)$), but
belonging to another class. Note that in the case of usual Jacobi
identity we have only cyclic permutations between three symbols
because the primary operations consist of two symbols and belong
to the unique class of cyclic permutations. The identity must have
symmetry with respect to the different classes. The identity is
expressed in the following form: \be

((A,B,C)D,E) + ((B,C,D)E,A) + ((C,D,E)A,B) + ((D,E,A)B,C) +
((E,A,B)C,D) + \\
+ ((C,B,A)E,D) + ((B,A,E)D,C) + ((A,E,D)C,B) + ((E,D,C)B,A) +
((D,C,B)A,E) + \\
+ ((D,A,C)E,B) + ((A,C,E)B,D) + ((C,E,B)D,A) + ((E,B,D)A,C) +
((B,D,A)C,E) + \\
+ ((C,A,D)B,E) + ((A,D,B)E,C) + ((D,B,E)C,A) + ((B,E,C)A,D) +
((E,C,A)D,B) + \\
 = 0

\ee For more detail about this identity see Appendix 2. Let us
formulate the following

\textbf{Proposition 2}: \textit{ For the 3-ary operation (12)
non-decomposable into binary ones there is a combinatorial
identity (18).}

This is similar to the existence for the matrix commutator Jacobi identity defined by the commutation operation.

\section{N-ary operation}

There are many various n-ary operations in the mixed tensor
algebra $\bigoplus\limits_{i,j=0}^\infty A_i^j$. For example,
there are 3-ary operation for $A^{2}_{1} \oplus A^{3}_{4}$ with
the five contractions, 5-ary for $A^{2}_{1} \oplus A^{2}_{5}$ with
the eight contractions and so on, and, possibly, the different
identities. The meaning of these identities as analogues of Jacobi identity
has the combinatorial nature. It pertains to the possibility to come to the same contraction diagram by different pathes having different weights. Under the "path" we mean the following: the word can be obtained by the different primary and secondary operations, and the sequence of this operations is called
the path.

Being restricted to the case of "linear" elementary operations for
some operator spaces let us consider the following table: \be

\begin{tabular}{|c|c|c|}
\hline
Operators&Space&Rank of the elementary operation\\
\hline
$A^{1}_{1}$ & $V^{1}$ & $2$ \\
\hline
$A^{2}_{1} \oplus A^{1}_{2}$ & $V^{1} \oplus V^{2}$ & $3$ \\
\hline
$A^{3}_{1} \oplus A^{1}_{2} \oplus A^{2}_{3}$ & $V^{1} \oplus V^{2} \oplus V^{3}$ & $4$ \\
\hline
$...$ & $...$ & $...$ \\
\hline $A^{n}_{1} \oplus A^{1}_{2} \oplus A^{2}_{3} \oplus ...
\oplus A^{n-1}_{n}$ &
$V^{1} \oplus V^{2} \oplus ... \oplus V^{n}$ & $n+1$ \\
\hline
\end{tabular}.

\ee In the column "Space" there are examples of the vector spaces
on which the tensors from column "Operator" act as operators.

Although the existence of the different generalized identities and
corresponding operations (composed by n-ary elementary operations)
is the subject of the further research let us define some basis
principles of construction the operations and identities. First of
all we give some remarks about the 3-ary operation (12). Defining
this operation we chose the second position in the sequence of the
objects $A^{2}_{1}$ and $A^{1}_{2}$ and assumed the following:
when $Y$ has the second position then the weight of the elementary
operation is $\alpha$, when $X$ - $\beta$, and when $Z$ -
$\gamma$. Note that we can chose the first and the third position.

For example, let us consider the operators of the space $A^{3}_{1}
\oplus A^{1}_{2} \oplus
 A^{2}_{3}$ and use the principle stated above.
In this case there is the primary operation consisted of four
symbols $(XYZW)$. After twice acting of this operation we get the
seven-symbolic words. Then the number of such operations is 35.
Choosing the first position we assume the following: when $X$ has
the first position then the weight of the elementary operation is
$\alpha$, when $Y$ - $\beta$, when $Z$ - $\gamma$ and when $W$ -
$\gamma$. We get that the number of the elementary operations with
the equal weight is 18 for the all type of the operators. In the
other hand the number of the elementary operations with the
different weight is 6 for the each type of the operators. The
total number of the elementary operations in each primary
operation is 72. In the identity the  number of terms must be $35
\cdot 2 = 72$. They are divided in groups having 7 terms in
according to the number of terms in each class of the cyclic
permutations. Each word must occur 8 times. We suppose that this
principles can be used in the cases of any n-ary operation.

\section{Discussion}

We suppose that the following questions must be interesting. First
of all there are general questions about operations. Given an
operation in ${\cal A}^\bullet_\bullet$ what is the set of
corresponding identities (in particular, for operations on ${\cal
A}^{1}_{2} \oplus {\cal A}^{2}_{1}$ having diagrams with
"non-linear" topology)? When can we restrict to double iterations
of the primary operation in the identity?

The operation (12) and the identity (18) allows to speak about the
existence of the structure on the space $A = A^{1}_{2} \oplus
A^{2}_{1}$ which is analogous to the Lie algebra on the matrix. So
there is the question about the analogues of Lie groups and
defining the exponential function in this case, if we start from
the notion of automorphism as the action which preserve the
structure. Note that 3-ary operation (12) and identity (18) can be
defined on the pairs of rectangle matrices. It is possible that
studying of this case is more simple then the case of tensors.

In the present paper we consider finite subalgebras of mixed
tensor algebra $\bigoplus\limits_{i,j=0}^\infty {\cal A}_i^j$. Note, that
in [4] there were considered infinite subalgebras and their
properties with respect to the nilpotent boundary operators
with $\delta^{k}=0, \ \ k \geq 3$. Another important class of
questions arises from the point of view when the mixed tensor
algebra is regarded as the one generated by a pair linear
spaces with bilinear form. This assumes the generalization in
which the free "algebras" are generated by $n$-tuples of
spaces contracted via $n$-linear forms.

\section{Appendix 1.}

Let us consider the operation (2) with values $\alpha = +1,\ \beta
= -1,\ \gamma = +1,\ \delta = -1$. We get: \be

A \circ B = A \cdot B - B \cdot A + A \cdot Tr B - B \cdot Tr A.

\ee Show that $(A \circ B) \circ C + (C \circ A) \circ B + (B
\circ C) \circ A \neq 0$. For this let us calculate: \be

(A \circ B) \circ C =\\
= + ABC - BAC - CAB + CBA - Tr C \cdot Tr A \cdot B + Tr C \cdot Tr B \cdot A - \\
- Tr A \cdot BC + Tr A \cdot CB + Tr B\cdot AC - Tr B \cdot CA
  + Tr C \cdot AB - Tr C \cdot BA \\

(C \circ A) \circ B =\\
= + CAB - ACB - BCA + BAC - Tr B \cdot Tr C \cdot A + Tr B \cdot Tr A \cdot C -\\
 - Tr C \cdot AB + Tr C \cdot BA + Tr A \cdot CB -
Tr A \cdot BC
  + Tr B \cdot CA - Tr B \cdot AC -\\
  - Tr B \cdot Tr C \cdot A + Tr B \cdot Tr A \cdot C\\

(B \circ C) \circ A =\\
= + BCA - CBA - ABC + ACB - Tr A \cdot Tr B \cdot C + Tr A \cdot Tr C \cdot B - \\
- Tr B \cdot CA + Tr B \cdot AC
  + Tr C \cdot BA - Tr C \cdot AB
  + Tr A \cdot BC - Tr A \cdot CB

\ee Summing we have the following expression: \be

(A \circ B) \circ C + (C \circ A) \circ B + (B \circ C) \circ A =\\

= Tr A \cdot (CB - BC)
 + Tr B \cdot (AC - CA)
  + Tr C \cdot (BA - AB).

\ee Now let us consider the result of the following expression
$((A \circ B) \circ C + (C \circ A) \circ B + (B \circ C) \circ A)
\circ D$: \be

((A \circ B) \circ C + (C \circ A) \circ B + (B \circ C) \circ A)
\circ D =\\
= - Tr A \cdot ( BC - CB) \cdot D
 + Tr B \cdot ( CA - AC ) \cdot D
  - Tr C \cdot ( AB - BA ) \cdot D -\\
  + D \cdot Tr A \cdot ( BC - CB)
 + D \cdot Tr B \cdot ( CA - AC )
  + D \cdot Tr C \cdot ( AB - BA ) -\\
- TrD \cdot Tr A \cdot ( BC - CB)
 - TrD \cdot Tr B \cdot ( CA - AC )
  - TrD \cdot Tr C \cdot ( AB - BA ).

\ee Summing this result with the results of the other expressions:
\be

((C \circ B) \circ D + ((D \circ C) \circ B + (B \circ D) \circ C)
\circ A,\\
((C \circ D) \circ A + (A \circ C) \circ D + (D \circ A) \circ C)
\circ B,\\
((B \circ A) \circ D + (D \cdot B) \circ A + (A \circ D) \circ B)
\circ C,

\ee
 we get the value equal to zero - (6).

\section{Appendix 2.}

Each term in (16) gives $36 \cdot 2$ words. The number of terms is
20, so the number of different words is $ \ $ - $720 \cdot 2$.
There are the words distinguished by the belonging to the
$(2,1)$-type or $(1,2)$-type tensors. So any word of each tensor
type occurs 6 times because the total number of different words is
$5! = 120$. Note that it is possible to compose 12 different words
using the different symbols, for example: $A^{2}_{1}, \ B^{1}_{2},
\ C^{2}_{1}, \ D^{1}_{2}, \ E^{2}_{1}$. It means that the set
having 120 different words divides in 10 classes. In each
five-symbolic word there are alternating tensors belonging to the
different types ((2,1) and (1,2)). So we can code the tensor type
by the tensor type of second and fourth symbol. At the other hand
there are 10 sets which consist of two symbol chosen from given
five ones. It defines the dividing of the set having 120 words of
each tensor type in 10 classes consisted of 12 words.

Now let us consider any class, for example $\{a,e\}$. We can form
the table in which to each word (left column) we associate the
primary operations (upper levels). It help us to find the term of
identity which participates in generating of considering word. In
lower levels there are the weights after two operations: \be

\begin{tabular}{|c|c|c|c|c|c|c|c|c|c|c|}
\hline
${A,E}$ & $abc,cba$ & $ace,eca$ & $cde,edc$ & $bda,adb$ & $aed,dea$ & $ebd,dbe$ & $dac,cad$ & $eab,bae$ & $ceb,bec$ & $ $\\
\hline
$baced$ & $+$ & $+$ & $+$ & $ $ & $ $ & $ $ & $ $ & $ $ & $ $ & $Eq1$\\
\hline
$badec$ & $ $ & $ $ & $+$ & $+$ & $+$ & $ $ & $ $ & $ $ & $ $ & $Eq2$\\
\hline
$bedac$ & $ $ & $ $ & $ $ & $ $ & $+$ & $+$ & $+$ & $ $ & $ $ & $Eq3$\\
\hline
$becad$ & $ $ & $+$ & $ $ & $ $ & $ $ & $ $ & $+$ & $ $ & $+$ & $Eq4$\\
\hline
$cabed$ & $+$ & $ $ & $ $ & $ $ & $ $ & $+$ & $ $ & $+$ & $ $ & $Eq2$\\
\hline
$cadeb$ & $ $ & $ $ & $ $ & $ $ & $+$ & $+$ & $+$ & $ $ & $ $ & $Eq3$\\
\hline
$cebad$ & $ $ & $ $ & $ $ & $+$ & $ $ & $ $ & $ $ & $+$ & $+$ & $Eq3$\\
\hline
$cedab$ & $ $ & $ $ & $+$ & $+$ & $+$ & $ $ & $ $ & $ $ & $ $ & $Eq2$\\
\hline
$dabec$ & $ $ & $ $ & $ $ & $+$ & $ $ & $ $ & $ $ & $+$ & $+$ & $Eq3$\\
\hline
$daceb$ & $ $ & $+$ & $ $ & $ $ & $ $ & $ $ & $+$ & $ $ & $+$ & $Eq4$\\
\hline
$decab$ & $+$ & $+$ & $+$ & $ $ & $ $ & $ $ & $ $ & $ $ & $ $ & $Eq1$\\
\hline
$debac$ & $+$ & $ $ & $ $ & $ $ & $ $ & $+$ & $ $ & $+$ & $ $ & $Eq2$\\
\hline
$к¬ми$ & $\beta \gamma$, $\gamma \alpha$ & $\alpha \beta$, $\alpha \beta$
& $\beta \gamma$, $\gamma \alpha$ & $\gamma \gamma$, $\alpha \beta$ & $\beta \beta$, $\beta \gamma$
& $\gamma \gamma$, $\alpha \beta$ & $\alpha \alpha$, $\alpha \gamma$ & $\beta \beta$, $\beta \gamma$
& $\alpha \alpha$, $\alpha \gamma$ & $ $\\
\hline
\end{tabular},

\ee where: \be

Eq1: 2 \cdot (\beta \gamma + \gamma \alpha + \alpha \beta) = 0\\
Eq2: \gamma \gamma + \gamma \alpha + \gamma \beta + \beta \beta +
\beta \alpha + \beta \gamma = 0\\
Eq3: \gamma \gamma + \beta \beta + \alpha \alpha + \gamma \alpha +
\gamma \beta + \beta \alpha = 0\\
Eq4: 2 \cdot (\alpha \alpha + \alpha \beta + \alpha \gamma) = 0.

\ee So, the total weight of each word (after summing) is equal to
zero. Note that this structure realizes in the other classes
parameterized by two symbols.

\paragraph{Acknowledgments}
Authors thank G.Sharigin for the useful remarks and A.Morozov for
fruitful discussions and attention to the work. The work of Yu.
Chernyakov was supported in part by grant RFFI ц 03-02-17554 and
grant ц 1999.2003.2 for support of scientific schools. The work of
V.Dolotin was supported by the grant RFFI ц 04-02-17227.

\end{document}